%

\input amstex
\documentstyle{amsppt}
%

%
%
\catcode`\@=11
\newif\ifdraft\draftfalse
\define\preliminary{\drafttrue}

%
%

\let\ga=\alpha

\let\gg=\gamma

\let\gk=\kappa
\let\gl=\lambda
\let\gs=\sigma
\let\gS=\Sigma

\let\gw=\omega


\define\ce{{\Cal E}}
\define\cf{{\Cal F}}
\define\cg{{\Cal G}}

\define\cm{{\Cal M}}

\define\catx#1#2{\mskip-2mu\raise#1\hbox{$#2\smallfrown$}\mskip-2mu}
\define\cat{\mathchoice{\catx{5pt}{\dsize}}{\catx{4.5pt}{\tsize}}
		{\catx{2.5pt}{\ssize}}{\catx{2pt}{\sssize}}}


\define\restrict{{\restriction}} 


\let\sat=\models

\define\sh{^{\mathchoice{}{\dsize\sharp}{\tsize\sharp}{\ssize\sharp}}}
 
\define\const#1{\frak c_{#1}}

\define\card#1{|#1|}
\define\Card#1{\Vert#1\Vert}
\define\set#1{\{\,#1\,\}}
\define\seq#1{(\,#1\,)}
\define\ve#1{\bold #1}
\define\image{\raise1.5pt\hbox{\lq\lq}\kern-.7pt} 
\predefine\integral{\int}


\define\ult{\operatorname{ult}}

\newcount\sectno
\global\sectno=0
\def\newsectno#1{\global\sectno=#1\global\thno=0\relax\the\sectno\ ---}
\outer\def\secthead#1#2{
	\global\thno=0\global\sectno=#1\subheading{\S. \the\sectno.\ #2}}

\def\mkinner#1#2{\edef#1{\expandafter\noexpand\csname#2\endcsname}}
\mkinner\innerproclaim{proclaim}
\mkinner\innerdefinition{definition}
\newcount\thno
\global\thno=0
\define\thethm{\ifnum\sectno=0\else\the\sectno.\fi\the\thno}
\define\doproclaim#1#2#3{\advance\thno by 1\relax
	#1{\ifNumberAfterThm#2 #3\else#3. #2\fi}}
\newif\ifNumberAfterThm
\NumberAfterThmtrue
\define\proc#1{\doproclaim{\innerproclaim}{#1}{\thethm}}

\outer\def\theorem{\proc{Theorem}}

\outer\def\prop{\proc{Proposition}}
\def\defin{\doproclaim{\innerdefinition}{Definition}{\thethm}}

\define\defining#1{{\it #1}}

\let\endth=\endproclaim


\def\putinmargin#1{\leavevmode\vadjust{
   \vbox to0pt{\vss\hbox to 0pt{\hss\hbox to0.70truein{\sevenrm#1\hss}}
   \vskip-8pt}}}
\newwrite\thmfile
\newif\iftagfile

\def\thmtag#1{\edef#1{\thethm}
	\ifdraft{\escapechar=-1\putinmargin{\string#1}}\fi}
\def\mytag#1#2{\expandafter\edef\csname#1\endcsname{#2}%
 \iftagfile\immediate\write\thmfile
        {\string\def\expandafter\string\csname#1\endcsname{#2}}\fi}
\def\rostertag#1{\mytag{#1}{\the\rostercount@}}
\def\nextsect#1/{{\def\xxx#1{\count0 by #1}\count0=\sectno
	\expandafter\advance\xxx#1\the\count0}}
\def\nextthm#1/{{\def\xxx#1{\count0 by #1}\count0=\thno
      \expandafter\advance\xxx#1\the\sectno.\the\count0}} 
\def\backthm#1{{\def\xxx#1{\count0 by -#1}\count0=\thno
	\expandafter\advance\xxx#1\relax\the\sectno.\the\count0}}
\define\proof{\demo{Proof}}
\define\endproofof#1{%
  \nobreak\hfil\penalty0\hfilneg\quad\qed\nobreak\quad\nobreak(#1)%
  \enddemo}
\define\rightqed{\hbox{}\nobreak\hskip0ptplus1filll\nobreak\hbox{}\nobreak\qed}
\define\endproof{\rightqed\enddemo}
\expandafter\ifx\csname today\endcsname\relax
	\def\today{\number\month/\number\day/\number\year}\fi
 
\def\lineto#1/{
	\def\xxx#1{height .4pt depth 0pt width #1}
	\expandafter\vrule\xxx#1
}



\catcode`\@=\active

\def\proc#1{
	\advance\thno by 1
	\innerproclaim{#1 \the\thno}
}
\def\defin{\advance\thno by1\innerdefinition{Definition \the\thno}}

\define\rank{\operatorname{rank}}
\define\concat#1#2{#1\cat#2}
\define\col{\Bbb P}
\define\sig#1#2{\gS^{#1}_{#2}}
\topmatter
\title The complexity of the core model\endtitle
\author William J. Mitchell\endauthor
\date\today\ifdraft--\preliminary\fi\enddate
\abstract We use the $\gS^1_3$~absoluteness  theorem to show that
complexity of the statement ``$(\gw,E)$ is isomorphic to an initial
segment of the core model'' is $\Pi^1_4$, and the complexity of the
statement ``$(\gw,E)$ is isomorphic to a member of $K$'' is $\Delta^1_5$.
\endabstract
\endtopmatter
\document

--- Complexity of ``$x$ \`is a real of $K$'' is $\Pi^1_3$. But we want
to look at sets of $K$ which are countable in $V$. 

\defin If $\Gamma$ is any point class over the reals then
well extend the usage of $\Gamma$ by saying that if $X\subset
H_{\gw_1}$ is in $\Gamma$ if it is each member of $X$ is transitive and 
$\set{E\subset{\gw^2}:(\exists x\in
X)(\gw,E)\cong (x,\in)}$ is in $\Gamma$.
\enddefinition

Suppose that $K=L[\ce]$, so that the initial segments of $K$ are the
sets of the form $J_{\ga}[\ce]$.

\prop
\roster
\item ``$\exists\ga\,``x=J_{\ga}[\ce]$'' is $\Pi^1_4$.
\item ``$x$ is in $K$'' is $\Delta^1_5$.
\endroster
\endth
\proof We consider clause~(1) first. The assertion that $x$ is a
premouse, that is, that there is an
ordinal $\ga$ and a good sequence $\cf$ of extenders such that 
$x=J_{\ga}(\cf)$, is $\gS^1_2$; and
the assertion that $J_{\ga}(\cf)$ is iterable is $\Pi^1_2$, so the
critical element in the characterization of initial segments of $K$ is
assertion that $\cf$ is maximal.  The assertion that 
$\cf$ is maximal for partial extenders, that is, that it contains
every mouse which it could have,
is ``every mouse $J_{\nu}[\cg]$ with projectum $\rho$ such that
$\cg\restrict\rho=\cf\restrict\rho$ is in $J_{\ga}[\cf]$,'' which is
$\Pi^1_3$ since
``$J_{\nu}[\cg]$ is a mouse'' is $\Pi^1_2$. 

Thus we need only consider the assertion that $\cf$ contains all of
the full measures which it should have.  By \cite{mit.coreii}, this 
means that $\cf$ contains every extender  $F$ satisfying
\roster
\item 
$F$ could be added to
$\cf\restrict\gg$ as $\ce_\gg$, that is, that $\concat{\cf\restrict\gg}
F$ is good at $\gg$, and
\item
If $J_{\nu}[\cg]$ is any iterable premouse such that
$\cg\restrict\gg=\ce\restrict\gg$ then $\ult(J_{\nu}[\cg],F)$ is well
founded.
\endroster
The assertion that an extender $F$ satisfies clause~(1) is $\gS^1_1$.
The assertion that a model $J_{\ga}[\cg]$ is iterable is $\Pi^1_2$, so
the assertion that an $F$ satisfies clause~(2) is $\Pi_{3}^1$ and
hence the assertion that every extender $F$ satisfying clauses~(1)
and~(2) is in $\cf$ is $\Pi^1_4$.
\medskip
Since a set $x\in H_{\gw_1}$ is a member of $K$ if and only if $x$ is
a member of some countable initial segment of $K$, it is easy express
``$x\in K$'' by a $\gS^1_5$ formula.  To see that ``$x\in K$'' can
also be expressed by a $\Pi^1_5$~formula, note that if $\nu$ is the
rank of $x$ (that is, the rank of the transitive closure of $x$
regarded as a well founded relation) then $x$ is in $K$ if and only if
it is in some mouse $J_{\ga}[\cg]$ such
$\cg\restrict\nu=\ce\restrict\nu$.   Thus $x\in K$ if and only if $$
\multline
\forall y\forall \nu\bigl((\nu=\rank(x)\text{ and }
y=\ce\restrict\nu)
\\
\qquad\implies\exists\cm\bigl(\cm=J_{\ga}[\cg] \text{ is an
iterable mouse with }\cg\restrict\nu=y\bigr)\bigr).
\endmultline\tag$*$ $$
By clause~(1) of the proposition, $y=\ce\restrict\nu$ is $\Pi^1_4$ and
hence ($*$) is $\sig15$.
\endproof

The main result of this paper is that the calculations in the last
proposition are best possible. We use the $\gS^1_3$-absoluteness
theorem from \cite{mit.abs}:

\theorem\thmtag\abs Suppose that there is no inner model of $\exists\gk
o(\gk)=\gk^{++}$ and that $a\sh$ exists for each real $a$.  Then any
model $M$ such that $K^{(M)}$ is an iterated ultrapower of $K$ is
$\gS^1_3$ correct.
\endth

\theorem\thmtag\mainthm\roster
\item
Suppose that $M$ is a model of set theory which satisfies that $\gk$
is measurable and that the sharp exists for every subset of $\gk^+$.
Then there is a model in which ``$x$ is an initial segment of $K$'' is
not expressible by a $\gS^1_4$ formula.
\item
Suppose that $M$ is a model of set theory which satisfies that there
are infinitely many measurable cardinals below $\gl$ and the sharp
exists for every subset of $\gl$.  Then there is a model in which
``$x\in K$'' cannot be expressed by a boolean combination of $\gS_4$
formulas.
\endroster\endth
\proof
Let $M$ satisfy the conditions of clause~(1). 
By taking a submodel if necessary we can assume that $M\sat\bigl(V=K$ and
$\lnot\exists\gk o(\gk)=\gk^{++}\bigr)$.  Let $\gk$ be the smallest
measurable cardinal in $M$, let $U$ be the measure on $\gk$ in $M$,
and let $i:M\to\ult(M,U)$ be the canonical embedding. Finally let
$\gl$ be the first fixed point of $i$ above $\gk$, so that
$\gk^+<\gl<\gk^{++}$, let $\col$ be the Levy collapse of $\gl$, and
let $G$ be $\col$-generic over $M$. We will show that there is no
$\gS^1_4$ formula $\phi$ such that
$$M[G]\sat\forall x\left(\phi(x)\iff \text{ $x$ is an initial segment
of $K$}\right).\tag$*$ $$

Suppose to the contrary that   $\phi$ is such a formula. Since $\phi$
doesn't contain any parameters, the homogeneity of the Levy collapse
implies that ($*$) is forced by the empty condition of $\col$.
Now $\col$ is also the Levy collapse of $\gl=i(\gl)$ in
$M_1=\ult(M,U)$, and $G$ is $\col$-generic over $M_1$, so ($*$) also
holds in $M_1[G]$. 

 Notice that $M_1[G]$ is a definable submodel of $M[G]$ and
$K^{(M_1)}=\ult(M,U)$ is an iterated ultrapower of $M=K^{(M)}$.
Furthermore, since $M$ satisfies that every subset of
$\gl$ has a sharp, it follows that in $M[G]$ every real has a
sharp. Thus
theorem~\abs, applied in $M[G]$, implies that every $\bold\gS^1_4$
formula true in $M_1$ is true in $M[G]$.  

Now let $x=J_{\gl}[\ce^{M_1}]$, so $x$ is
countable in the models $M$ and $M_1$.  
Since $x$ is an initial segment of $K^{(M_1[G])}$ but not of
$K^{(M[G])}$ we know that 
$M_1[G]\sat\phi(x)$, while $M\sat\lnot\phi(x)$, but since $\phi$ is
$\gS^1_4$ this contradicts the last paragraph, and this contradiction
completes the proof of clause~(1).

\medskip
In order to prove clause~(2) suppose that $M$ is a model of set theory
containing infinitely many measurable cardinals, and
that if $\gl$ is the sup of the first $\gw$~measurable cardinals in
$M$ then $x\sh\in M$ for each $x\subset\gl$ in $M$.
As before, we can assume that $M\sat
\bigl(V=K+\lnot\exists\gk\,o(\gk)=\gk^{++}\bigr)$, so that $M=L[\ce]$.
Let $\col$ 
be the partial order to collapse $\gl$, let $G$ be $\col$-generic over
$M$, and let $a$ be a real which is Cohen generic over $M[G]$. Our
model will be $M_0=M[G,a]$. We will show first that $x\in K$ is not
expressible in $M_0$ by any formula of the form $\pi(x)\land\gs(x)$
where $\pi(x)$ is a $\Pi^1_4$ formula and $\gs$ is a
$\gS^1_4$-formula, and afterward we will generalize this to arbitrary
boolean combinations of $\gS^1_4$ formulas.
 	
\define\ceonegl{\ce_1\restrict\gl}

Suppose then, that $x\in K$ is expressible in $M_0$ by such a formula
$\pi(x)\land\gs(x)$. 
Let $(U_i:i\in\gw)$ be the first $\gw$ measures in $M=L[\ce]$, let
$i_1:L[\ce]\to L[\ce_1]$ be the iterated ultrapower using each of the
measures $U_i$ exactly once, let $i_a:L[\ce]\to L[\ce_a]$ be the
iterated ultrapower using only the measures $U_i$ for $i\in a$, and
Set $N=L[\ce_a][G,a]$ and $M_1=L[\ce_1][G,a]$. 
Then $$M_0\supset N\supset M_1.$$ 
Then $\ce_1\restrict\gl$ is a member of $K^{(M_0)}$ and $K^{(M_1)}$,
but not a member of $K^{(N)}$.  Just as in the proof of clause~(1),
the formula ($*{}*$) holds in $N$ and
$M_1$, so $$\pi(\ce_1\restrict\gl)\land\gs(\ce_1\restrict\gl)\tag$*{*}*$ $$ is true
in $M_0$ and $M_1$ but false in $N$.  But 
lemma~\abs\ implies that $M_1$ is $\sig13$~correct in $N$, so that 
$$M_1\sat\gs(\ceonegl)\quad\implies\quad N\sat\gs(\ceonegl)$$
and $N$ is $\sig13$~correct in $N$ so 
$$M_1\sat\pi(\ceonegl)\quad\implies\quad N\sat\pi(\ceonegl),$$
contradicting formula~($*{*}*$) and hence disproving $(*{*})$.
\medskip
\define\ceigl{\ce_i\restrict\gl}
Now suppose that  $$M_0\sat\forall x(x\in K\iff\phi(x))$$ where
$\phi(x)$ is a boolean combination of $\sig14$ formulas. Then
$\phi(x)$ is equivalent to a formula of the form
$$\bigvee_{i<n}(\pi_i(x)\land\gs_i(x)$$ for some $n\in\gw$ and
formulas $\pi_i(x)$ and $\gs_i(x)$ which are $\Pi^1_3$ and $\sig13$,
respecitively. Now we define a chain $$M_0\subset N_0\subset
M_1\subset\dots\subset N_{n-1}\subset M_n$$ of models. As before, we
have $M_0=M[G,a]=L[\ce][G,a]$.  Let $j_i:L[\ce]\to L[\ce_i]$ be the
iterated ultrapower using each of the measures $\seq{U_i:i\in\gw}$
exactly $i$ times, and let $L[\ce_{i,a}]$ be the result of iterating
each of the measures $\set{j_{i}(U_k):k\in a}$ once more. As before,
$\ce_n$ will be a member of $M_i$ for each $i\le n$, and formula~(4)
is true for each of the models $M_i$ so that $\phi(\ceigl)$ is true in
each $M_i$. Then by the pigeon hole principle there are integers $k<n$
and $i_0<i_1\le n$   such that $\gs_k(\ceigl)\land\pi_k(\ceigl)$ is
true in both the models $M_{i_{0}}$ and $M_{i_1}$. Just as in the last
paragraph this implies that $\gs_k(\ceigl)\land\pi_k(\ceigl)$ is true
in $N_{i_0}$, but this is impossible because $\ce_{n}$ is not in
$K^{(N_{i_0})}$.
\endproof

The hypothesis cannot be weakened in theorem~\mainthm. Taking
clause~(1) for example, if $N$ is any model such that $K^{(N)}$ doesn't
have any measurable cardinals which are countable in $N$ then $x$ is
an initial segment of $K^{(N)}$ if it contains all the mice it should,
which is a
$\Pi^1_3$ statement; while if $K^{(N)}=L[\ce]$ where $\ce$ is countable in $N$ then
$\ce$ is characterized by the fact that $\ce\sh$ doesn't exist and
$\ult(L_{\gw_1}(\ce),E)$ is ill founded for any ultrafilter $E$ on
$L[\ce]$ which is not either on the sequence $\ce$ or an iterate of a measure
which is on $\ce$.  The hypothesis to clause~(1) is the weakest
which will allow the existence of a model $N$ such that $K^{(N)}$
contains a countable mouse and is not equal $L[\ce]$ for a sequence
$\ce$ which is countable in $N$.

Similarly, the boldface version of theorem~\mainthm\ requires that the
measurable cardinals of $K^{(N)}$ be cofinal in $\gw_1^{(N)}$, which
requires that we start with an inaccessible limit of measurable
cardinals.  The arguments above readily yield
\theorem Suppose that there is a model  $M$ of set theory with an an
inaccessible limit of measurable cardinals. Then there is a model $N$
such that  ``$x$ is an initial segment of $K$'' is
not expressible in $N$ by a $\boldsymbol\gS^1_4$ formula, and ``$x\in K$''
cannot be expressed in $N$  by any boolean combination of $\boldsymbol\gS_4$
formulas.
\endth

\enddocument